\documentclass[12pt,reqno]{amsart}

\usepackage{amssymb,latexsym}
\usepackage[matrix,arrow]{xy}

\newtheorem{thm}{Theorem}[section]

\newtheorem{lem}[thm]{Lemma}

\theoremstyle{definition}
\newtheorem{defn}{Definition}[section]

\theoremstyle{remark}
\newtheorem{rem}{Remark}[section]
\newcommand{\n}{\nabla}
\newcommand{\ns}{\nabla^{\mathcal{S}}}
\newcommand{\np}{\nabla^{\partial M}}
\newcommand{\e}{\eqno}
\newcommand{\vi}{\varphi}
\newcommand{\s}{\Bbb{S}}
\newcommand{\ms}{\mathcal{S}}
\newcommand{\md}{\mathcal{D}}
\newcommand{\gp}{\gamma^{\partial M}}
\newcommand{\gs}{\gamma^{\mathcal{S}}}
\newcommand{\p}{\Bbb{S}|_{\partial M}}
\newcommand{\ld}{\lambda}

\newcommand{\g}{\gamma}
\newcommand{\G}{\Gamma}

\begin{document}
\title{\text Eigenvalue estimates  for Dirac operator with the generalized APS boundary condition \bf  }
\author []{ Daguang Chen }
\email { chendg@amss.ac.cn } \curraddr{Institute of Mathematics,
 Academy of Mathematics and Systems Sciences,
 Chinese Academy of Sciences, Beijing 100080, P.R. China.}
\thanks{}
\date{21 Feb-2006}
\begin{abstract}
Under two boundary conditions: the generalized
\-Atiyah-\-Patodi-Singer boundary condition and the modified
generalized \-Atiyah-\-Patodi-Singer boundary condition, we get
the lower bounds for the eigenvalues of the fundamental Dirac
operator on compact spin manifolds with nonempty boundary.
\end{abstract}
\maketitle
\renewcommand{\sectionmark}[1]{}
\section{introduction}
Let $M$ be a compact \-Riemannian spin manifold without boundary.
In 1963, \-Lichnerowicz \cite{L} proved that any eigenvalue of the
Dirac operator satisfies
$$\ld^2>\frac 14 \underset{M}\inf R,$$
where $R$ is the positive scalar curvature of $M$. The first sharp
estimate for the smallest  eigenvalue $\ld$ of the Dirac operator
$D$ was obtained by \-Friedrich \cite{Fr2} in 1980. The idea of
the proof is based on using a suitably modified \-Riemannian spin
connection. He proved the inequality
$$
\ld^2\geq \frac{n}{4(n-1)}\underset{M} \inf R,
$$
on closed manifolds $(M^n, g)$ with the positive scalar curvature
$R$. Equality  gives an Einstein metric. In 1986, combining the
technique of the modified spin connection with a conformal change
of the metric, \-Hijazi \cite{H1} showed, for $n\geq3$
$$
\ld^2\geq\frac {n}{4(n-1)} \mu_1,
$$
where $\mu_1$ is the first eigenvalue of the conformal \-Laplacian
given by $L=\frac{4(n-1)}{n-2}\triangle+R$. When the equality
holds, the manifold becomes an Einstein manifold.

Let $M$ be a compact \-Riemannian spin manifold with nonempty
boundary. In this case, the boundary conditions for \-spinors
become crucial to make the Dirac operator elliptic, and, in
general, two types of the boundary conditions are considered, i.e.
the local boundary condition and the \-Atiyah-\-Patodi-Singer
(\-APS) boundary condition (see, e.g.\cite{FS}). In 2001,
\-Hijazi, \-Montiel and \-Zhang \cite{HMZ1} proved the generalized
version of \-Friedrich-type inequalities under both the local
boundary condition and the \-Atiyah-\-Patodi-Singer boundary
condition. In 2002, \-Hijazi, \-Montiel and Rold\'{a}n \cite{HMR}
investigated the \-Friedrich-type inequality for eigenvalues of
the fundamental Dirac operator under four elliptic boundary
conditions and under some curvature assumptions (the non-negative
mean curvature). In particular, they introduced a new global
boundary condition, namely the modified \-Atiyah-\-Patodi-Singer
(\-mAPS) boundary condition.

In the present paper, we study the spectrum of the fundamental
Dirac operator on compact \-Riemannian spin manifolds with
nonempty boundary, under the two new boundary conditions: the
generalized \-Atiyah-\-Patodi-Singer (\-gAPS) boundary condition
and the modified generalized \-Atiyah-\-Patodi-Singer (\-mgAPS)
boundary condition. Following the terminology of $[6]$ those are
obtained composing a so-called spectral projection with the
identity in the first case, and with the zero order differential
operator $\mathrm{Id}+\gamma(e_0)$ in the second case (where
$\gamma(e_0)$ is the Clifford product with the unit normal $e_0$
on the boundary $\partial M$).

\section{Preliminaries}

Let $(M,g)$ be an (n+1)-\-dimentional \-Riemannian spin manifold
with nonempty boundary $\partial M$. We denote by $\n$ the
Levi-Civit\`{a} connection on the tangent bundle $TM$. For the
given structure $Spin M$ (and so a corresponding orientation) on
manifold $M$, we denote by $\s$ the associated \-spinor bundle,
which is a complex vector bundle of rank $2^{[\frac {n+1}2]}$.
Then let
$$
\gamma :\Bbb {C}l(M)\longrightarrow End_{\Bbb C}(\Bbb {S})
$$
be the Clifford multiplication, which is a fibre preserving
algebra morphism. It is well known \cite{LM} that there exists a
natural Hermitian metric $(,)$ on the '\-spinor bundle $\Bbb S$
which satisfies the relation
$$
(\gamma (X)  \varphi,\gamma (X)  \psi)=|X|^2(\varphi,\psi),\e(2.1)
$$
for any vector field $X \in \Gamma(TM)$, and for any \-spinor
fields $\varphi,\psi\in \Gamma(\Bbb S)$. We denote also by $\n$
the \-spinorial Levi-Civit\`{a} connection acting on the \-spinor
bundle $\s$. Then the connection $\n$ is compatible with the
Hermitian metric $(,)$ and Clifford multiplication "$\g$".
Recalling the fundamental Dirac operator $D$ is the first order
elliptic differential operator acting on the \-spinor bundle $\s$,
which is locally given by
$$
D=\sum_{i=0}^{n}\g (e_i)\n_i,\e(2.2)
$$
where $\{e_0,e_1,\cdots,e_n\}$ is a local \-orthonormal frame of
$TM$.

Since  the normal bundle to the boundary \-hypersurface is
trivial, the \-Riemannian manifold $\partial M$ is also spin and
so we have a corresponding \-spinor bundle  $\Bbb S\partial M$, a
\-spinorial Levi-Civit\`{a} connection $\np$ and an intrinsic
Dirac operator $D^{\partial M}=\gp(e_i)\np_i$. Then a simple
calculation yields the \-spinorial Gauss formula
$$
\n_i\vi=\np_i\vi+\frac12 h_{ij}\g(e_j)\g(e_0)\vi,
$$
for $1\leq i,j\leq n, \vi\in \G(\s|_{\partial M})$, $e_0$ is a
unit normal vector field compatible with the induced orientation,
$h_{ij}$ is the second fundamental form of the boundary
\-hypersurface.

It is easy to check (see \cite{B}\cite{HM}\cite{HMR}\cite{HMZ2})
that the restriction of the \-spinor bundle $\s$ of $M$ to its
boundary is related to the intrinsic Hermitian \-spinor bundle
$\Bbb S\partial M$  by
$$
\ms:=\p\cong
\begin{cases}\Bbb S\partial M &\mbox{if $n$ is even}\\
\Bbb S\partial M\oplus\Bbb S\partial M&\mbox{if $n$ is odd}.
\end{cases}
$$
For any \-spinor field $\psi\in\G({\ms})$ on the boundary
$\partial M$, define on the restricted \-spinor bundle $\ms$, the
Clifford multiplication $\gs$ and the \-spinorial connection $\ns$
by
$$
\gs(e_i)\psi:=\g(e_i)\g(e_0)\psi\e(2.3)
$$
$$
\ns_i\psi:=\n_i\psi-\frac 12 h_{ij}\gs(e_j)\psi=\n_i\psi-\frac
12h_{ij} \g(e_j)\g (e_0)\psi.\e(2.4)
$$
One can easily find that $\ns$ is compatible with the Clifford
multiplication $\gs $, the induced Hermitian inner product $( , )$
from $M$ and together with the following additional identity
$$
\ns_i(\g(e_0)\psi)=\g(e_0)\ns_i\psi.\e(2.5)
$$
As a consequence, the boundary Dirac operator $\md $ associated
with the connection $\ns$ and the Clifford multiplication $\gs$,
is locally given by
$$
\md\psi=\gs(e_i)\ns_{e_i}\psi=\g(e_i)\g(e_0)\n_{e_i}\psi+\frac H
2\psi,\e(2.6)
$$
where $H=\sum h_{ii}$ is the mean curvature of $\partial M$. In
fact, $\g(e_0)(\mathcal D -\frac12H)=\g(e_i)\n_i $ is the
\-hypersurface Dirac operator defined by \-Witten \cite{W} to
prove the positive energy conjecture in general relativity.

Now by (2.5), we have the \-supersymmetry property
$$
\md\g(e_0)=-\g(e_0)\md.\e(2.7)
$$
Hence, when $\partial M$ is compact, the spectrum of $\md$ is
symmetric with respect to zero and coincides with the spectrum of
$D^{\partial M}$ for $n$ even and with (Spec$D^{\partial
M}$)$\cup$(-Spec$D^{\partial M}$) for $ n$ odd.

From \cite{HM} or \cite{HMR}, one gets the integral form of the
Schr\"{o}\-dinger-\-Lichnerowicz formula for a compact
\-Riemannian spin manifold with compact boundary
$$
\aligned
 \int_{\partial M}(\varphi,\mathcal{D}\varphi)-\frac 1 2 \int_{\partial
M}H|\varphi|^2=\int_M|\n\varphi|^2+\frac R 4
|\varphi|^2-|D\varphi|^2.
\endaligned
\eqno(2.8)
$$

By $(2.7)$,  $\mathcal D$ is self-adjoint with respect to the
induced Hermitian metric $(,)$ on $\ms$. Therefore, $\mathcal{D}$
has a discrete  spectrum contained in $\Bbb R$ numbered like
$$
\cdots\leq\ld_{-j}\leq\cdots\ld_{-1}<0\leq\ld_0\leq
\ld_1\leq\cdots\leq\ld_j\leq\cdots
$$
and one can find an \-orthonormal basis $\{\vi_j\}_{j\in\Bbb Z}$
of $L^2(\partial M;\ms)$ consisting of  eigenfunctions of $\md$
(i.e.\,$\md\vi_j=\ld_j\vi_j, j\in \Bbb Z$) (see e.g. \textbf{Lemma
1.6.3} in \cite{G}). Such a system $\{\ld_j;\vi_j\}_{j\in\Bbb Z}$
is called a $spectral$ $decomposition$ of $L^2(\partial M;\ms)$
generated by $\md$, or, in short, a $spectral$ $resolution$ of
$\md$. According to the definition $14.1$ in \cite{BW}, we have
 \begin{defn} For self-adjoint Dirac operator $\mathcal
 D$ and for any real $b$ we shall denote by
 $$
P_{\geq b}:L^2(\partial M;\ms)\longrightarrow L^2(\partial M;\ms)
 $$
the $spectral\, projection$, that is, the orthogonal projection of
$L^2(\partial M;\ms)$ onto the subspace spanned by \{
$\varphi_j|\ld_j\geq b$\}, where {$\{\ld_j;\varphi_j\}_{j\in
\Bbb{Z}}$ } is a spectral resolution of $\mathcal D$.
\end{defn}
\begin{defn}
For $\vi\in\Gamma(\mathcal S)$ and $b\leq 0$, the projective
boundary conditions
$$
\begin{aligned}
&P_{\geq b}\vi=0, \\
&P^m_{\geq b}\vi=P_{\geq b}(Id+\g(e_0))\vi=0
\end{aligned}
$$
are called the generalized Atiyah\--Patodi\--Singer (\-gAPS)
boundary condition and the modified generalized
\-Atiyah-\-Patodi-Singer (\-mgAPS) boundary condition
respectively.
\end{defn}
\begin{rem}
 We shall adopt the notation $P_{<b}=Id-P_{\geq
b}$ ; $P_{[b,-b]}=P_{\geq b}-P_{>-b}$ for $b<0$ or
$P_{[-b,b]}=P_{\geq -b}-P_{>b}$ for $ b>0$. It is well known that
the spectral projection $P_{\geq b}$ is a pseudo-differential
operator of order zero (see e.g. the proposition 14.2 in \cite
{BW}).
\end{rem}
\begin{rem} If $b=0$, the \-gAPS
 boundary condition is exactly the
APS boundary condition (see \-Atiyah,\-Patodi
 \& Singer\cite{APS1}\cite{APS2}\cite{APS3}).
\end{rem}
 In case of a closed manifold, the fundamental Dirac operator is a formally
 self-adjoint operator and so its spectrum is discrete and real. While in case
 of a manifold with nonempty boundary, we have the following
 formula
$$
\int_M (D\varphi,\psi)-\int_M(\varphi,D \psi)=-\int_{\partial
M}(\g(e_0)\varphi,\psi).\eqno(2.9)
$$
for $\vi,\psi\in\G(\partial M)$, and where $e_0$ is the inner unit
field along the boundary. In the following, we'll show ellipticity
and self-\-adjointness for the Dirac operator under both the
\-gAPS boundary condition and the \-mgAPS boundary condition.

On one hand, if $P_{\geq b}\varphi =0,P_{\geq b}\psi =0$, for
$\varphi,\psi\in\Gamma(\ms)$ and $b\leq 0$, then
$P_{<b}(\g(e_0)\varphi)=\g(e_0) P_{>-b} $. Therefore
$$
\int_{\partial M}(\g(e_0)\varphi,\psi)=0.\e(2.10)
$$
On the other hand, for $\vi,\psi\in\Gamma(\ms)$ and $b\leq 0$,
under the \-mgAPS boundary
 conditions, i.e.
 $$P_{\geq b}(Id+\g(e_0))\vi=0\qquad \mbox{and} \qquad P_{\geq
 b}(Id+\g(e_0))\psi=0,$$
a simple calculation implies that (2.10) also holds.

 These facts imply that the fundamental Dirac operator $D$ is self-adjoint
under the \-gAPS  boundary condition and the \-gmAPS boundary
condition and hence it has
real and discrete eigenvalues.\\

The ellipticity of both \-gAPS and \-mgAPS boundary conditions the
fundamental Dirac operator $D$ follows straightforward from $[6,
\textrm{Prop. }14.2]$, $[16,\textrm{Prop. }1]$ and $[16,
\textrm{Thm. }5]$.\\

\begin{lem}
Under the \emph{\-gAPS} boundary condition  $P_{\geq b}\vi=0$, the
inequality
$$
\int_{\partial M}(\vi,\md\vi)\leq b\int_{\partial M}|\vi|^2
$$
holds for $b\leq 0$.
\end{lem}
\begin{proof}
 Let $\{\ld_j;\varphi_j\}_{j\in
\Bbb{Z}}$ be a spectral resolution of the \-hypersurface Dirac
operator $\mathcal D$. Then any $\vi\in\G(\ms)$ can be expressed
as follows,
$$
\vi=\textstyle\sum_j c_j\vi_j,
$$
where $c_j=\int_{\partial M}(\vi,\vi_j)$. Then we have
$$
\int_{\partial M} (\varphi,\mathcal
D\varphi)=\underset{\ld_j<b}{\textstyle\sum}\ld_j|c_j|^2\leq
b\underset{\ld_j<b}{\textstyle\sum}|c_j|^2=b\int_{\partial
M}|\vi|^2.
$$
\end{proof}
\section{lower bounds for the eigenvalues}

In this section, we adapt the arguments used in \cite{HZ1} and
\cite{HMZ1} to the case of compact \-Riemannian spin manifolds
with nonempty boundary. We use the integral identity (2.8)
together with an appropriate modification of the Levi-Civit\`{a}
connection to obtain the eigenvalue estimates.
\begin{thm}
Let $M^n$ be a compact \-Riemannian spin manifold of dimension
$n\geq 2$, with the nonempty boundary $\partial M$, and let $\ld$
be any eigenvalue of $D$ under the \emph{\-gAPS} boundary
condition $P_{\geq b}\varphi=0$, for $\varphi \in  \Gamma(\ms)$
and $b\leq0$. If there exist real functions $a,u$ on $M$ such that
$$
b+ a\ du(e_0)\leq \frac{1}{2}H
$$
on $\partial M$, then
$$
\ld^2>\frac{n}{4(n-1)}\underset{a,u}\sup \ \underset{M}\inf
R_{a,u},\e(3.1)
$$
where
$$R_{a,u}:=R-4a\Delta u+4\n a\cdot\n
u-4(1-\frac{1}{n})a^2|du|^2\e(3.2)
$$
$\Delta$ is the positive scalar \-Laplacian, R is scalar curvature
of $M$.
\end{thm}
\begin{proof}
For any real function $a$ and $u$, define the modified
 \-spinorial connection (see \cite{HMZ1}\cite{HZ1}) on $\G{(\Bbb S)}$ by
 $$
\nabla_i^{a,u}=\n_i+a\n_iu+\frac{a}{n}\n_ju\g(e_i)
\g(e_j)+\frac{\ld}{n}\g(e_i).\eqno(3.3)
 $$
 A simple calculation yields
$$
\aligned
\int_M|\nabla^{a,u}\varphi|^2&=\int_M\Big[(1-\frac{1}{n})\ld^2
-\frac{R_{a,u}}{4}\Big]|\varphi|^2\\
&+\int_{\partial M}(\varphi,\mathcal{D}\vi)+\int_{\partial
M}\Big[a\ du(e_0)-\frac{H}{2}\Big]|\varphi|^2.
\endaligned\eqno(3.4)
$$
Considering \textbf{Lemma 2.1} and $ b+ a\ du(e_0)\leq
\frac{1}{2}H$, we infer
$$
\ld^2\geq\frac{n}{4(n-1)}\underset{a,u}\sup \ \underset{M}\inf
R_{a,u}.\e(3.5)
$$
If the equality in (3.5) holds, then by the \textbf{Lemma 3} in
\cite{HMZ1}, we have $a=0$ or $u=constant$ and $b\leq \frac12 H$.
Since $\vi$ is a non-degenerated killing \-spinor, $|\varphi|^2$
is a nonzero constant. Let $\{\ld_j;\varphi_j\}_{j\in \Bbb{Z}}$ be
a spectral resolution of the \-hypersurface Dirac operator
$\mathcal D$. Under the \-gAPS boundary condition $P_{\geq
b}\vi=0$, one gets $\vi=\underset{\ld_j<b}{\textstyle\sum}
c_j\vi_j $, where $c_j=\int_{\partial M}(\vi,\vi_j)$. Then we have
$$
\begin{aligned}
0=&\int_{\partial M}(\vi,\md \vi)-\int_{\partial M}\frac
12H|\vi|^2\\
=&\underset{\ld_j<b}{\textstyle\sum}\ld_j |c_j|^2 -\frac
12\underset{\ld_j,\ld_k<b}\sum c_jc_k\int_{\partial M}
H(\vi_j,\vi_k)\\
\leq&\underset{\ld_j<b}{\textstyle\sum}(\ld_j-b)|c_j|^2<0.
\end{aligned}
$$
This is a contradiction. Thus the inequality (3.1)
holds.
\end{proof}
\begin{rem}If $b=0$, the \-gAPS boundary condition becomes
the \-APS boundary condition and then the theorem is exactly the
\textbf{Theorem 3.1} in \cite{HMZ1}.
\end{rem}


Now we make use of the energy-momentum tensor, introduced in
\cite{H2} and used in \cite{HMZ2}\cite{HZ1}\cite{HZ2}, to get
lower bounds for the eigenvalues of $D$. For any \-spinor field
$\varphi\in\G(\s)$, we define the associated energy-momentum
2-tensor $Q_{\varphi} $ on the complement of its zero set by
$$
Q_{\varphi,ij}=\frac{1}{2}\Re\big(\g(e_i)\n_j\varphi
+\g(e_j)\n_i\varphi,\varphi/|\varphi|^2\big).
$$
Obviously, $Q_{\varphi,ij}$ is a symmetric tensor. If $\varphi$ is
the \-eigenspinor of the Dirac operator $D$, the tensor
$Q_\varphi$ is well-defined in the sense of distribution.

\begin{thm} Let $M^n$ be a compact \-Riemannian spin manifold of
dimension $n\geq 2$, whose boundary $\partial M$ is nonempty, and
let $\ld$ be any eigenvalue of $D$ under the \emph{\-gAPS}
boundary condition $P_{\geq b}\varphi=0,$ for $b\leq 0,\varphi\in
\Gamma(\ms)$. If there exist real functions $a,u$ on $M$ such that
$$
b+a du(e_0)\leq H/2
$$
on $\partial M$, where $H$ is the mean curvature of $\partial M$,
then
$$
\ld^2\geq \underset{a,u}\sup\ \underset{M}\inf(\frac
14R_{a,u}+|Q_\varphi|^2),
$$
where $R_{a,u}$ is given in (3.2).
\end{thm}
\begin{proof}
For any real function $a$ and $u$, define the modified \-spinorial
connection (see \cite{HMZ1}\cite{HZ1}) on $\G(\s)$ by
$$
\nabla_i^{Q,a,u}=\n_i+a\n_i u+\frac{a}{n}\n_ju\g(e_i)\g(e_j)
+Q_{\varphi,ij}\g(e_j).\e(3.6)
$$
One can easily compute
$$
\aligned \int_M|\nabla^{Q,a,u}\varphi|^2
=&\int_M\Big[\ld^2-(\frac{R_{a,u}}{4}+|Q_\varphi|^2)
\Big]|\varphi|^2\\
&+\int_{\partial M}(\varphi,\mathcal{D}\varphi)+\int_{\partial
M}\Big[a du(e_0)-H/2\Big]|\varphi|^2.
\endaligned\e(3.7)
$$
Considering the boundary conditions, we immediately arrive at the
asserted formula.
\end{proof}
\begin{rem}
By \-Cauchy-\-Schwarz inequality and $tr Q_\vi=\ld$, we have
$$
|Q_\vi|^2=\underset{i,j}{\textstyle \sum} |Q_{\vi,ij}|^2\geq
\underset{i}{\textstyle \sum} |Q_{\vi,ii}|^2\geq \frac 1n(tr
Q_\vi)^2=\frac 1n \ld^2.
$$
Therefore, one gets the inequality (3.5).
\end{rem}
\begin{rem}If $b=0$, the above theorem becomes \textbf{Theorem 5}
 in \cite{HMZ1}. Under the \-gAPS boundary condition,
 taking $a=0$ or $u=constant$ in (3.6) and
assuming $H\geq 0$, then one gets the following inequality
\cite{H2}
$$\ld^2\geq \underset{M}\inf\big( \frac R4+|Q_\vi|^2\big).$$
\end{rem}
In the following, we state the eigenvalue estimates for the Dirac
operator under the \-mgAPS boundary condition.

\begin{thm}
Let $M^n$ be a compact \-Riemannian spin manifold of dimension
$n\geq 2$, whose boundary $\partial M$ is nonempty, and let $\ld$
be any eigenvalue of $D$ under the \emph{\-mgAPS} boundary
condition $P_{\geq b}^m=P_{\geq b}(Id+\g(e_0))=0,$ for $b\leq 0$.
If there exist real functions $a,u$ on $M$ such that
$$
a\ du(e_0)\leq \frac{1}{2}H   \e(3.8)
$$
on $\partial M$, then
$$
\ld^2\geq\frac{n}{4(n-1)}\underset{a,u}\sup \ \underset{M}\inf
R_{a,u},\e(3.9)
$$
where $R_{a,u}$ is given in (3.2).
 Moreover, the equality holds
if and only if $M$ carries a nontrivial real killing \-spinor
field with a killing constant $-\frac\ld n<\frac b{n-1}$ and the
boundary $\partial M$ is minimal.
\end{thm}
\begin{proof} Assuming
$$P^m_{\geq b}\varphi=P_{\geq b}(\varphi+\g(e_0)\varphi)=0\qquad
\mbox{and}\qquad P^m_{\geq b}\psi=P_{\geq b}(\psi+\g(e_0)\psi)=0,
$$
for $\varphi,\psi \in \Gamma (\ms)$ and $b\leq0$, then we have
$$
P_{>-b}\psi+\g(e_0)P_{<b}\psi=0\quad \mbox{and}\quad
P_{[b,-b]}\psi+\g(e_0)P_{[b,-b]}\psi=0,
$$
i.e.
$$
P_{\geq b}^m\psi=P_{>-b}\psi+\g(e_0)P_{<b}\psi=0\quad
\mbox{and}\quad P_{[b,-b]}\psi=0 \eqno(3.10)
$$
These imply
$$
\begin{aligned}
P^m_{\geq b}\varphi=P_{>-b}(\varphi+\g(e_0)\varphi)=0,\\
P^m_{\geq b}\psi=P_{>-b}(\psi+\g(e_0)\psi)=0.
\end{aligned}
$$
Then one gets
$$
\g(e_0)\psi-\psi =P_{>-b}(\g(e_0)\psi-\psi). \eqno(3.11)
$$
The relation (2.7) implies that
$$
(\mathcal{D}\psi,\psi)=\frac {1}{2}(\mathcal{D}(\psi+\g(e_0)
\psi),\psi-\g(e_0) \psi).
$$
The combination $P^m_{\geq
b}\psi=P^m_{>-b}\psi=P_{>-b}(\psi+\g(e_0)\psi)=0$ with (3.11)
yields
$$
\int_{\partial M}(\mathcal{D}(\psi+\g(e_0)\psi),\psi-\g(e_0)
\psi)=0.
$$
This implies
$$
\int_{\partial M} (\mathcal{D}\psi,\psi)=0.\eqno(3.12)
$$
From (3.3), (3.4) and (3.12), one can infer
$$
\int_M|\nabla^{a,u}\varphi|^2=\int_M\Big[(1-\frac{1}{n})\ld^2
-\frac{R_{a,u}}{4}\Big]|\varphi|^2 +\int_{\partial M}\Big[a\
du(e_0)-\frac{H}{2}\Big]|\varphi|^2. \eqno(3.13)
$$
Then the inequality (3.9) holds.

If the equality occurs, we deduce
$$\n^{a,u}\vi=0\qquad  \mbox{and}\qquad \frac{1}{2}H=a du(e_0).$$
By the \textbf{Lemma 3} in \cite{HMZ1}, we get $a=0$ or
$u=constant$. Therefore
$$H=0\qquad \mbox{and}\qquad\n_i\vi=-\frac{\ld}{n}\g(e_i)\vi.$$
From the \-supersymmetry property (2.7), we get
$$
\mathcal{D}(\varphi+\g(e_0)\varphi)=-\frac {n-1}n \ld
(\varphi+\g(e_0) \varphi).
$$
Since $P^m_{\geq b}\varphi=0$, we deduce that
$-\frac{n-1}{n}\ld<b$.

Conversely, if $M$ is a compact \-Riemannian spin manifold with
minimal boundary $\partial M$ and a nontrivial killing \-spinor
$\varphi$ with a real killing constant $-\ld/n<\frac{b}{n-1}$,
then we have $D\vi=\ld\vi$. Moreover, from the fact that
$\n_{e_0}\varphi=-\frac\ld n\g(e_0)\varphi$, we infer that the
restriction of $\varphi$ to the boundary satisfies
$$
\mathcal{D}\varphi=-\frac{n-1}{n}\g(e_0)\varphi.
$$
Finally, we have
$$
\mathcal{D}(\varphi+\g(e_0)\varphi)=-\frac{n-1}{n}\ld(\varphi+\g(e_0)\varphi).
$$
Since  the \-spinor field $\varphi+\g(e_0)\varphi$ is an
\-eigenspinor of $\mathcal{D}$ with a eigenvalue
$-\frac{(n-1)}n\ld<b$, one can infer $P^m_{\geq b}\varphi=0$.
\end{proof}
\begin{rem} Under the \-mgAPS boundary condition, taking $a=0$ or
$u=constant$ in (3.8) and (3.13), one gets \-Friedrich's
inequality
$$
\ld^2\geq \frac{n}{4(n-1)}\underset{M}\inf R.
$$
In particular, if we take $b=0$, the above theorem is exactly the
\textbf{Theorem 5} in \cite{HMR}.
\end{rem}

Using the energy-momentum tensor, (3.7) and (3.12), we get
$$
 \int_M|\nabla^{Q,a,u}\varphi|^2
=\int_M\Big[\ld^2-(\frac{R_{a,u}}{4}+|Q_\varphi|^2)
\Big]|\varphi|^2+\int_{\partial M}\Big[a
du(e_0)-H/2\Big]|\varphi|^2.
$$
Thus we obtain the following theorem:
\begin{thm}
Let $M^n$ be a compact \-Riemannian spin manifold of dimension
$n\geq 2$, with nonempty boundary $\partial M$, and let $\ld$ be
any eigenvalue of $D$ under the \emph{\-mgAPS} boundary condition
$P_{\geq b}^m\varphi=P_{\geq b}(\varphi+\g(e_0)\varphi)=0,$ for
$b\leq 0$, $\varphi\in \Gamma(\ms)$. If there exist real functions
$a,u$ on $M$ such that
$$
a\ du(e_0)\leq \frac{1}{2}H
$$
on $\partial M$, then
$$
\ld^2\geq\underset{a,u}\sup \ \underset{M}\inf \big(\frac
14R_{a,u}+|Q_\vi|^2\big),\e(3.14)
$$
where $R_{a,u}$ is given in (3.2).
\end{thm}

\textbf{Acknowledgments} The author would like to  express his
gratitude to Professor X. \-Zhang for his suggestions and support.
He would also like to thank the referee for many valuable
suggestions.


\providecommand{\bysame}{\leavevmode\hbox
to3em{\hrulefill}\thinspace}

\end{document}